\magnification1200

\input amssym

\font\tenrm=pplr9d\font\sevenrm=pplr9d at 7.6pt\font\fiverm=pplr9d at 6pt\rm

\font\mathtext=pplri9o \fontdimen2\mathtext=0pt
 \font\mathsubtext=pplri9o at 7.6pt\fontdimen2\mathsubtext=0pt
 \font\mathsubsubtext=pplri9o at 6pt\fontdimen2\mathsubsubtext=0pt
\font\mathlet=eurmo10\fontdimen2\mathlet=0pt
 \font\mathsublet=eurmo10 at 7.6pt  \fontdimen2\mathsublet=0pt
 \font\mathsubsublet=eurmo10 at 6pt \fontdimen2\mathsubsublet=0pt
\skewchar\mathlet='177 \skewchar\mathsublet='177 \skewchar\mathsubsublet='177
\textfont1=\mathlet
 \scriptfont1=\mathsublet
 \scriptscriptfont1=\mathsubsublet
\textfont0=\mathtext
 \scriptfont0=\mathsubtext
 \scriptscriptfont0=\mathsubsubtext
\textfont4=\tenrm
 \scriptfont4=\sevenrm
 \scriptscriptfont4=\fiverm

\mathcode`\(="4428
\mathcode`\)="5429
\mathcode`\:="343A
\mathcode`\;="643B
\mathcode`\[="445B
\mathcode`\]="545D
\mathcode`\+="242B
\mathchardef\colon="643A

\mathchardef\Gamma="0100
\mathchardef\Delta="0101
\mathchardef\Theta="0102
\mathchardef\Lambda="0103
\mathchardef\Xi="0104
\mathchardef\Pi="0105
\mathchardef\Sigma="0106
\mathchardef\Upsilon="0107
\mathchardef\Phi="0108
\mathchardef\Psi="0109
\mathchardef\Omega="010A

\font\tenbbm=bbmsl10
\font\sevenbbm=bbmsl10 at 8.6pt
\font\fivebbm=bbmsl10 at 6pt
\newfam\bbmfam
\textfont\bbmfam=\tenbbm
\scriptfont\bbmfam=\sevenbbm
\scriptscriptfont\bbmfam=\fivebbm

\font\tt=txtt
\font\sc=pplrc9d

\parindent=0pt
\parskip=\smallskipamount

\def\frac#1/#2{\leavevmode\kern.1em
\raise.5ex\hbox{\the\scriptfont0 #1}\kern-.1em
/\kern-.15em\lower.25ex\hbox{\the\scriptfont0 #2}}

\def\proclaim #1. #2\par{\medbreak
  \noindent{\sc#1.\enspace}{\rm#2}\par\medbreak}
\def\gh{generalized homology}

\centerline{\bf Odd-dimensional Charney-Davis Conjecture}
\bigskip
\centerline{\char'221wiatos\char'252aw
R.~Gal\footnote{$^a$}{Partially supported by
Polish {\sc n201 012 32/0718} grant.}
\& Tadeusz Januszkiewicz\footnote{$^b$}{Partially supported by the
{\sc nsf} grant {\sc dms-0706259}.}}
\footnote{}{2000 {\it Mathematics Subject Classification: }\sc 52b70
(52b11, 06a07).}
\footnote{}{{\it Key phrases: flag complex, h-vector,
Charney--Davis Conjecture}.}
\bigskip

{
{\sc Abstract:}
More than once we have heard that the Charney-Davis Conjecture
makes sense only for odd-dimensional spheres. This is to point out
that in fact it is also a statement about even-dimensional spheres.}
\bigskip

A conjecture of Heinz Hopf asserts that
the sign of the Euler characteristic of a smooth Riemannian $2d$-dimensional
manifold of non-positive sectional curvature is the same for all such
manifolds, this is the same as that of product of
non-positively curved surfaces:
$$(-1)^d\chi (M^{2d})\geq 0.$$

For Riemannian manifolds the condition of non-positive sectional
curvature is equivalent to being locally {\sc cat(0)}.
The Hopf Conjecture subsequently has been generalized to include
closed, piecewise Euclidean locally, {\sc cat(0)} (\gh) manifolds.

\bigskip
By work of M.~W.~Davis [D], Coxeter groups provide a rich source of piecewise
Euclidean, locally {\sc cat(0)} spaces.
Given a {\it flag} triangulation of a (\gh) sphere $L^{n-1}$,
a construction of Davis gives a reflection (\gh) orbifold
${\cal O}^n$, with many (generalized homology) manifold covers.

The Euler characteristic of $\cal O$ is given by the f- and h-polynomials
of $L$ as follows:
$$\chi({\cal O})=f_L(-1/2)=h_L(-1).$$

Charney and Davis emphasized the combinatorial implications of
the Hopf Conjecture in this context.

\proclaim Conjecture ({\rm [CD, Conj.~D, p.~135]}).
If $L$ is a flag triangulation of the (\gh) sphere
of dimension $2d-1$ then $(-1)^dh(-1)\geq 0$.

\bigskip
The Hopf Conjecture does not say anything about odd-dimensional manifolds.
So on the face of it, the Charney-Davis Conjecture should not say anything
about odd-dimensional (\gh) spheres.
The point we want to make in this note is that in fact it does.

\proclaim Theorem. The Charney-Davis Conjecture is equivalent to
the following statement.
Let $L$ be a {\gh} sphere of dimension $2d$.
Let $h_L(t)$ be its h-polynomial, and let $\widetilde h_L(t)$ be defined by
$h_L(t)=(1+t)\widetilde h_L(t)$. Then
$$(-1)^d\widetilde h_L(-1)\geq 0.$$

{\sl Proof:}
Let $L$ be a suspension of $\char'212$.
Since the h-polynomial is multiplicative for joins
$(1+t)\widetilde h_L(t)=h_L(t)=(1+t)h_{\char'212}(t)$.
The statement $(-1)^dh_{\char'212}(-1)=(-1)^d\widetilde h_L(-1)\geq 0$ is just
the Charney-Davis Conjecture for $\char'212$.

To prove the other implication
recall that f-polynomial and h-polynomial are related by the formula
$$(2+t)(1+t)^{2d-1}\widetilde h_L\left({1\over1+t}\right)=
(1+t)^{2d}h_L\left({1\over1+t}\right)=t^{2d}f_L\left({1\over t}\right).$$
Differentiating both sides we get
$$(2+t)\left[(1+t)^{2d-1}\widetilde h_L\left({1\over1+t}\right)\right]'
+(1+t)^{2d-1}\widetilde h_L\left({1\over1+t}\right)
=2d\,t^{2d-1}f_L\left({1\over t}\right)-t^{2d-2}f'\left({1\over t}\right)$$
where $\left[(1+t)^{2d-1}\widetilde h_L(1/(1+t))\right]'$ is a polynomial.
Substitute $t=-2$ and use the fact that, by Dehn-Sommerville,
$f_L(-1/2)=0$ to get
$$(-1)^{2d-1}\widetilde h_L(-1)=-(-2)^{2d-2}f'(-1/2).$$

We omit the proof of the following straightforward claim.
The sum of f-polynomials of links of vertices of $L$
is equal to the derivative of the f-polynomial of $L$.

Applying the above claim to the preceding equality gives
$$(-1)^d\widetilde h_L(-1)=4^{d-1}\sum_v (-1)^d h_{\mathop{\rm Lk}_v}(-1).$$
The right hand side is non-negative by the Charney-Davis Conjecture.
Hence the proof.
\hfill$\square$

\proclaim Remark. The quantity $(-1)^d\widetilde h_L(-1)$ is equal
to $\gamma_d(L)$, the top coefficient of the
$\gamma$-polynomial introduced in [G, Def.~2.1.4]. The calculation
proving that $\gamma_d(L)>0$ provided the Charney-Davis.~Conjecture holds
for links of all vertices in $L$ was mentioned in [G, Cor.~2.2.2]
without relating it to the h-polynomial of $L$.

In view of the above the Charney-Davis Conjecture for even-dimensional
spheres is essentially included in (though perhaps a dramatic restatement of)
the Conjecture 2.1.7 in [G] which treats equally even- and odd-dimensional
spheres and provides further strengthenings of the Charney-Davis Conjecture.

One may speculate about the geometric interpretation of $\widetilde h_L(-1)$.
Note that presumably $\widetilde h_L$ is an h-polynomial of a
$(2d-1)$-dimensional sphere. For example, if $L$ is a icosahedron,
$\widetilde h_L$ is an h-polynomial of the decagon. On the other hand
the geometry of the Davis orbifolds for the icosahedron
and a suspension of the decagon are very different. The former
is hyperbolic and the latter is a product. Thus there is no hope of
relating $\widetilde h_L(-1)$ to $\ell^2$-torsion.

\bigskip
\centerline{\sc References}\smallskip

\parindent=.4in

\item{[CD]} {\sc R.~Charney \& M.~Davis}, {\it The Euler characteristic
of a non-positively curved, piecewise Euclidean manifold},
Pacific J.~Math. {\bf 171} (1995), pp.~117--137,

\item{[D]} {\sc Michael~W.~Davis}, {\it The geometry and topology of Coxeter
groups}, {\it London Mathematical Society Monographs Series}, {\bf 32}. Princeton University Press, Princeton, {\sc nj}, 2008,

\item{[G]} {\sc \char'221.~R.~Gal}, {\it Real Root Conjecture fails for five
and higher dimensional spheres}, Discrete \& Computational Geometry {\bf 34}
(2005), pp.~269--284.

\bigskip
\obeylines\parskip=0pt
{\sc \char'221wiatos\char'252aw R.~Gal:}
Mathematical Institute, Wroc\char'252aw University
pl. Grunwaldzki \frac2/4, 50-384 Wroc\char'252aw, Poland
{\tt sgal@math.uni.wroc.pl}
\medskip

{\sc Tadeusz Januszkiewicz:}
Department of Mathematics, The Ohio State University
231 {\sc w 18th Ave}, Columbus, {\sc oh 43210, usa}
and the Mathematical Institute of Polish Academy of Sciences;
on leave from Mathematical Institute, Wroc\char'252aw University
{\tt tjan@math.ohio-state.edu}

\bigskip

\rightline{\it Columbus, 12. May 2009}
\bye